\documentclass[11pt,hyp,reqno]{amsart}
\usepackage{amsmath,amsfonts,a4,epsfig}
\usepackage{graphicx}
\usepackage[latin1]{inputenc}
\usepackage[french, english]{babel}
\usepackage{color}
\usepackage{enumerate}
\usepackage{setspace}
\usepackage{fancyhdr}

\def\preuve{\begin{proof}}

\newtheorem{defi}{Definition}[section]

\newtheorem{lemm}{Lemma}[section]
\newtheorem{prop}{Proposition}[section]
\newtheorem{rem}{Remark}[section]

\newtheorem{comm}{Comment}[section]
\newtheorem{coro}{Corollary}[section]
\newtheorem{theo}{Theorem}[section]

  \newenvironment{demo}{\noindent {\it Proof:}
      \begin{quotation}\noindent}{\end{quotation}\hfill$\square $}
\usepackage{hyperref}
\hypersetup{nesting=true,debug=true,naturalnames=true}
\usepackage{graphicx,amssymb,upref}

\let\<\langle
\let\>\rangle
\usepackage[all,pdf]{xy}
\UseComputerModernTips

\let\uml\"

\title{Non self-adjoint Laplacians on a directed graph}

\author{ MARWA BALTI}
\address{Universit\'e de Carthage, Facult\'e des Sciences de Bizerte: Math\'ematiques et Applications (UR/13ES47) 7021-Bizerte (Tunisie)\\
Universit\'e de Nantes, Laboratoire de Math\'ematique Jean Lauray, CNRS, Facult\'e des Sciences, BP 92208, 44322 Nantes, (France).
}
\email{ balti-marwa@hotmail.fr}



\begin{document}
\begin{abstract}
We consider a non self-adjoint Laplacian on a directed graph with non symmetric edge weights. We analyse spectral properties of this Laplacian under a Kirchhoff assumption. Moreover we establish isoperimetric inequalities in terms of the numerical range to show the absence of the essential spectrum of the Laplacian on \textit{heavy end} directed graphs. %We introduce a special self-adjoint operator and compare its essential spectrum with that of the non self-adjoint Laplacian considered.
\end{abstract}
\maketitle
\tableofcontents

\section*{Introduction}
 The non self-adjoint operators are more difficult to study than the self-adjoint ones:  no spectral theorem in general, wild resolvent growth...  The related  theory is studied by different authors: L. N. Trefethen \cite{Tr} for non symmetric matrices, W. D. Evans, R. T. Lewis, A. Zettl \cite{Ev} and R. T. Lewis \cite{Lws} for non self-adjoint operators in a Hilbert space. Recently, the interest in spectral properties of non self-adjoint operators has already led to a variety of new results, both in the continuous and discrete settings, e.g, bounds on complex eigenvalues \cite{R.F} and Lieb-Thirring type inequalities \cite{Han}, \cite{DHK}.
  This can be explained by the complicated structure of the resolvent of such an operator seen as an analytic function. In this paper we focus on directed graphs to study a non symmetric Laplacian. We develop a general approximation theory for the eigenvalues on directed graphs with non symmetric edge weights assuming only a condition of "total conductivity of the vertices" presented as the Assumption $(\beta)$. We investigate the spectrum of our discrete non self-adjoint Laplacian. We collect some basic properties of the Laplacian and we seek to show the emptiness of its essential spectrum by using isoperimetric inequalities. We explain how isoperimetric inequalities  can be linked to the numerical range of non symmetric operators.
   In fact, for the self-adjoint Laplace-Beltrami operator, Jeff Cheeger proved an inequality that links the first nontrivial eigenvalue on a compact Riemannian manifold to a geometric constant $h$. This inspired an analogous theory for graphs (see \cite{Fuji}, \cite{A.G}). In this work, we introduce a kind of Cheeger constant on a filtration of a directed graph $G$ and we estimate the associated Laplacian $\Delta$. We give an estimation for the numerical range of $\Delta$ in terms of the Cheeger constant. We use this estimation and propose a condition on the weights for the absence of essential spectrum of \textit{heavy end} directed graphs. There is an analogous result of H. Donnelly and P. Li \cite{Don} for a self-adjoint operator on complete negatively curved manifolds. They show that the Laplacian on a rapidly curving manifold has a compact resolvent.\\
Section \ref{In} is devoted to some definitions and notions on a directed graph with non symmetric edge weights and the associate non symmetric differential Laplacian. We describe some basic results: Green's formula and the spectral properties of $\Delta$ and of its formal adjoint.\\
In Section \ref{2}, we study spectral properties of the bounded operator $\tilde{\Delta}$ by relying on known results for the symmetric case. \\
 In Section \ref{ch}, we establish the Cheeger inequality for the non symmetric Dirichlet Laplacian on any subset of the set of vertices $V$ to give a lower bound for the bottom of the real part of the numerical range. We control the real part of the numerical range of $\Delta$ and relate it with the spectrum of the its closure $\overline{\Delta}$. We characterize the absence of essential spectrum of $\overline{\Delta}$. Fujiwara \cite{Fuji} and Keller \cite{Klr} introduced a criterion for the absence of essential spectrum of the symmetric Laplacian on a rapidly branching graph. In fact, our criterion is: positivity of the Cheeger constant at infinity on \textit{heavy end} graphs.
\section{Preliminaries}\label{In}
We review in this section some basic definitions on infinite weighted graphs and introduce the notation used in the article. They are introduced in \cite{Ma} for finite non symmetric graphs (see  \cite{AT} and \cite{N.T} for the symmetric case).
\subsection{Notion of Graphs}
A directed weighted graph is a triplet $G :=(V,\vec{E},b)$, where $V$ is a countable set (the vertices), $\vec{E}$ is the set of directed edges and $b : V \times V\to\left[0,\infty\right)  $ is a weight function satisfying the following conditions:
\begin{itemize}
  \item $b(x,x)=0$ for all $x\in V$
  \item $b(x,y)>0$ iff $(x,y)\in \vec{E}$
    \end{itemize}

  In addition, we consider a measure on V given by a positive function
$$m:V\to (0,\infty) .$$
The weighted graph is \textit{symmetric} if  for all $x,y\in V$, $b(x,y)=b(y,x)$,
  as a consequence $(x,y)\in \vec{E}\Leftrightarrow (y,x)\in\vec{E}$.\\
  The graph is called  \textit{simple} if the weights $m$ and $b$ are constant and equal to $1$ on $V$ and $\vec{E}$ respectively.\\

The set $E$ of undirected edges is given by
$$E=\left\{\{x,y\},~(x,y)\in\vec{E}\text{ or } (y,x)\in\vec{E}\right\}.$$%, \text{ and } x\sim y \text{ if } \{x,y\}$$
\begin{defi}
Define for a subset $\Omega$ of $V$, the vertex boundary and the edge boundary of $\Omega$ respectively by:
$$\partial_V\Omega=\big\{ y\in\Omega:~\{x,y\}\in E~\text{ for some } x\in \Omega^c\big\} $$
$$\partial_E \Omega=\Big\{(x,y)\in \vec{E}:~ (x\in \Omega,~y\in \Omega^{c} )~~or ~~(x\in \Omega^{c},~y\in \Omega)\Big\}.$$
\end{defi}

On a non symmetric graph we have two notions of connectedness.
\begin{defi}
\begin{itemize}
\item A path between two vertices $x$ and $y$ in $V$ is a finite set of directed edges $(x_1,y_1);~(x_2,y_2);..;(x_n,y_n),~n\geq 2$ such that
$$x_1=x,~y_n=y \text{ and } x_i=y_{i-1}~~\forall~2\leq i\leq n$$
\item $G$ is called connected if two vertices are always related by a path.
\item $G$ is called strongly connected if for all vertices $x,y$ there is a path from $x$ to $y$ and one from $y$ to $x$.
 \end{itemize}
\end{defi}

\subsection{Functional spaces}
Let us introduce the following spaces associated to the graph $G$:\\
\begin{itemize}

\item the space of functions on the graph $G$ is considered as the space of complex
functions on V and is denoted by
$$\mathcal{C}(V)=\{f:V\to \mathbb{C} \}$$
\item $\mathcal{C}_c(V)$ is its subset of finite supported functions; \\
 
\item we consider for a measure $m$, the space
$$\ell^2(V,m)=\{f\in \mathcal{C}(V), ~~\sum_{x\in V}m(x)|f(x)|^2<\infty\}.$$
It is a Hilbert space when equipped by the scalar product given by
$$(f,g)_m=\sum_{x\in V}m(x)f(x)\overline{g(x)}.$$
The associated norm is given by:
$$\|f\|_m=\sqrt{(f,f)_m}.$$
%\item  For $U\subset V$, we define
 %$$\mathcal{C}_c(U)=\{f:U\to \mathbb{C}, \text{ with finite support}\}.$$
 
\end{itemize}
\subsection{Laplacian on a directed graph}
In this work, we assume that the graph under consideration is connected, locally finite, without loops and satisfies for all $x\in V$ the following conditions:
$$\sum_{y\in V}b(x,y)>0 \text{ and }\sum_{y\in V}b(y,x)>0.$$
 We introduce the combinatorial Laplacian $\Delta$
defined on $\mathcal{C}_c(V)$ by:
$$\Delta f(x)=\frac{1}{m(x)}\sum_{y\in V}b(x,y)\left( f(x)-f(y)\right) .$$
For all $x\in V$ we note by $\beta^{+}(x)=\displaystyle{\sum_{y\in V}b(x,y)}$, in particular if $m(x)=\beta^+(x)$ then the Laplacian is said to be the normalized Laplacian and it is defined by:
 $$\tilde{\Delta}f(x)=\frac{1}{\beta^+(x)}\sum_{y\in V}b(x,y)\big(f(x)-f(y)\big).$$

   \textbf{ Dirichlet operator}: Let $U $ be a subset of $ V$, $f\in \mathcal{C}_c(U)$ and $g:V\to \mathbb{C}$ the extension of $f$ to $V$ by setting $g = 0$ outside $U$.
 For any operator $A$ on $\mathcal{C}_c(V)$, the Dirichlet operator $A^{D}_U$ is defined by 
   $$ A^{D}_U(f)=A(g)\vert_ U .$$
   
The operator $\Delta$ may be non symmetric if the edge weight is not symmetric.
 \begin{prop}
The formal adjoint $\Delta'$ of the operator $\Delta$ is defined on $\mathcal{C}_c(V)$ by:
$$\Delta' f(x)=\frac{1}{m(x)}\left( \sum_{y\in V}b(x,y)f(x)-\sum_{y\in V}b(y,x)f(y)\right) .$$
\end{prop}
\begin{demo}
For all $f,g \in \mathcal{C}_c(V)$, we have
 \begin{align*}
(\Delta f,g)_m=&\sum_{(x,y)\in \vec{E}}b(x,y)\big( f(x)-f(y)\big) \overline{g(x)}\\
=&\sum_{x\in V}f(x)\sum_{y\in V}b(x,y)\overline{g(x)}-\sum_{(y,x)\in \vec{E}}b(y,x) \overline{g(y)}f(x)\\
=&\sum_{x\in V}f(x) \left( \sum_{y\in V}b(x,y)\overline{g(x)}-\sum_{(y,x)\in \vec{E}}b(y,x) \overline{g(y)}\right).
\end{align*}
As $(\Delta f,g)_m=(f,\Delta'g)_m$, so we get
$$\Delta' f(x)=\frac{1}{m(x)}\left( \sum_{y\in V}b(x,y)f(x)-\sum_{y\in V}b(y,x)f(y)\right) .$$
\end{demo}
\begin{rem}
The operator $\Delta'$ can be expressed as a Schr\"odinger operator with the potential $q(x)=\displaystyle{\frac{1}{m(x)}\sum_{ y\in V} \big(b(x, y)-b(y, x)\big),~ x\in V}$:
$$\Delta' f(x)=\frac{1}{m(x)}\sum_{y\in V}b(y,x)\big(f(x)-f(y)\big)+q(x)f(x).$$
\end{rem}

We introduce here the \textbf{Assumption} $(\beta)$ and we assume that it is satisfied by the considered weighted graph, throughout the rest.\\

  \textbf{Assumption $(\beta)$}: for all $x\in V$,~ $\beta^+(x)=\beta^-(x)$\\
  where $$\beta^+(x)=\sum_{y\in V }b(x,y)\text{ and }\beta^-(x)=\sum_{y\in V}b(y,x).$$
    \begin{rem}
The Assumption $(\beta)$ is natural. Indeed, it looks like the Kirchhoff's law in the electrical networks.
\end{rem}
\begin{coro}\label{comm} 
We suppose that the Assumption $(\beta)$ is satisfied, the operator $\Delta'$  is simply a Laplacian, given by

$$\Delta' f(x)=\frac{1}{m(x)}\sum_{y\in V}b(y,x)\big(f(x)-f(y)\big).$$
\end{coro}
In the sequel, for the sake of simplicity we introduce the symmetric Laplacian $H$ associated to the graph
with the symmetric edge weight function $a(x,y) = b(x,y)+b(y,x)$. It acts on $\mathcal{C}_c(V)$ by,
$$Hf(x)=(\Delta+\Delta')f(x)=\frac{1}{m(x)}\sum_{y\in V}a(x,y)\big(f(x)-f(y)\big).$$ 
The quadratic form $Q_\Delta$ of $H$ is given by
$$Q_\Delta(f)=(\Delta f,f)+\overline{(\Delta f,f)},~f\in \mathcal{C}_c(V).$$

\begin{comm}
Let $f\in \mathcal{C}_c(V)$, we have $Q_\Delta(f)=2\mathcal{R}e(\Delta f,f)$. Then
\begin{equation}\label{symm}
\inf_{\|f\|_m=1}Q_\Delta(f)=\inf_{\|f\|_m=1}2\mathcal{R}e(\Delta f,f).
\end{equation}
\end{comm}
We establish an explicit Green's formula associated to the non self-adjoint Laplacian $\Delta$ from which any estimates on the symmetric quadratic form $Q_\Delta$ can be directly cited from the literature.
\begin{lemm}(Green's Formula)\label{grn}
Let $f$ and $g$ be two functions of $\mathcal{C}_c(V)$. Then under the Assumption $(\beta)$ we have
$$(\Delta f,g)_m+\overline{(\Delta g,f)}_m=\sum_{(x,y)\in \vec{E}}b(x,y)\big( f(x)-f(y)\big) \big( \overline{g(x)-g(y)}\big).$$
\end{lemm}
\begin{demo}
The proof is given by a simple calculation. From Corollary \ref{comm}, we have
\begin{align*}
(\Delta f,g)_m+\overline{(\Delta g,f)}_m=&(H f,g)_m\\
=&\sum_{(x,y)\in \vec{E}} b(x,y)\big(f(x)-f(y)\big)\overline{g(x)}\\
+&\sum_{(y,x)\in \vec{E}} b(y,x)\big(f(x)-f(y)\big)\overline{g(x)}\\
=&\sum_{(x,y)\in \vec{E}} b(x,y)\Big(f(x)\overline{g(x)}+f(x)\overline{g(x)}-f(y)\overline{g(x)} -f(x)\overline{g(y)}\Big)\\
=&\sum_{(x,y)\in \vec{E}} b(x,y)\big( f(x)-f(y)\big) \big(\overline{g(x)-g(y)}\big).
\end{align*}
\end{demo}

We refer to \cite{Kts} page 243 for the definitions of the spectrum and the essential spectrum of a closed operator $A$ in a Hilbert space $\mathcal{H}$, with domain $D(A)$.
\begin{defi}\label{deff}
\begin{itemize}
\item The spectrum $\sigma(A)$ of $A$ is the set of all complex numbers $\lambda$ such that $(A-\lambda)$ has no bounded inverse.
\item The essential spectrum $\sigma_{ess}(A)$ of $A$ is the set of all complex numbers $\lambda$ for which the range $R(A-\lambda)$ is not closed or $\dim\ker(A-\lambda)=\infty$.
\end{itemize}
\end{defi}
\section{Spectral analysis of the bounded case}\label{2}
This part concerns some basic properties of the bounded non self-adjoint Laplacian $\tilde{\Delta}$. We introduce the concept of the numerical range. It has been extensively studied the last few decades. This is because it is very useful in studying and understanding the spectra of operators (see \cite{Ber}, \cite{JY}, \cite{Arl}).\\
%In the mathematical field of linear algebra and convex analysis, it is important to introduce the notion of the numerical range.
\begin{defi}
 The numerical range of an operator $T$ with domain $D(T)$, denoted by  $W(T)$ is the non-empty set
$$W(T)=\{(T f,f),~~f\in D(T),~\parallel f\parallel=1 \}.$$
\end{defi}

   The  following  Theorem in \cite{JY} shows  that  the  spectrum behave nicely with respect to the closure of the numerical range.

\begin{theo}
Let $\mathcal{H}$ be a reflexive Banach space and $T$ a bounded operator on $\mathcal{H}$. Then: $$\sigma(T)\subset \overline{W(T)}.$$
\end{theo}
The following Proposition is one of the main tools when working with the normalized Laplacian.
\begin{prop}
Suppose that the Assumption $(\beta)$ is satisfied. Then $\tilde{\Delta}$ is bounded by 2.
\end{prop}
\begin{demo}
\begin{align}\label{C.S}
\mid(\tilde{\Delta} f,g)_{\beta^+}\mid=&\mid \sum_{x\in V}\overline{g(x)}\sum_{y\in V} b(x,y)\big(f(x)-f(y)\big)\mid\notag\\
\leq& \sum_{x\in V}\beta^+(x) |f(x)\overline{g(x)}|+\sum_{x\in V}\mid\overline{g(x)}\mid\sum_{y\in V} b(x,y)\mid f(y)\mid\notag\\
\end{align}
by the Assumption $(\beta)$ and the Cauchy-Schwarz inequality we prove the result:
\begin{align*}
\mid(\tilde{\Delta} f,g)_{\beta^+}\mid\leq& \big(f,f\big)_{\beta^+}^{\frac{1}{2}}\big(g,g\big)_{\beta^+}^{\frac{1}{2}}+\sum_{x\in V}\mid\overline{g(x)}\mid\Big(\sum_{y\in V} b(x,y)\Big)^{\frac{1}{2}}\Big(\sum_{y\in V} b(x,y)\mid f(y)\mid^2\Big)^{\frac{1}{2}}\\
\leq& \big(f,f\big)_{\beta^+}^{\frac{1}{2}}(g,g)_{\beta^+}^{\frac{1}{2}}+\Big(\sum_{x\in V}\beta^+(x)\mid g(x)\mid^2\Big)^{\frac{1}{2}}\Big(\sum_{x\in V}\sum_{y\in V} b(x,y)\mid f(y)\mid^2\Big)^{\frac{1}{2}}\\
\leq& (f,f)_{\beta^+}^{\frac{1}{2}}(g,g)_{\beta^+}^{\frac{1}{2}}+(g,g)_{\beta^+}^{\frac{1}{2}}\Big(\sum_{y\in V}\mid f(y)\mid^2\sum_{x\in V} b(x,y)\Big)^{\frac{1}{2}}\notag\\
\leq& \big(f,f\big)_{\beta^+}^{\frac{1}{2}}\big(g,g\big)_{\beta^+}^{\frac{1}{2}}+\big(g,g\big)_{\beta^+}^{\frac{1}{2}}\Big(\sum_{y\in V}\mid f(y)\mid^2\beta^-(y)\Big)^{\frac{1}{2}}\notag\\
\leq& 2\big(f,f\big)_{\beta^+}^{\frac{1}{2}}\big(g,g\big)_{\beta^+}^{\frac{1}{2}}\notag.
\end{align*}
Then $$\|\tilde{\Delta}\|_{\beta^+}=\sup_{\|f\|_{\beta^+}\leq 1\atop \|g\|_{\beta^+}\leq 1}\mid(\tilde{\Delta} f,g)_{\beta^+}\mid\leq 2.$$
\end{demo}

It is useful to develop some basic properties of the numerical range to make the computations of the spectrum of the Laplacian.
\begin{prop}
\begin{enumerate} Let $G$ be a connected graph, satisfying the Assumption $(\beta)$. Then
  \item $\sigma(\tilde{\Delta})\subset  D(1,1)$, the closed disc with center $(1,0)$ and radius $1$.
  \item If $\beta^+(V)< \infty $, then $0$ is a simple eigenvalue of $\tilde{\Delta}$.
\end{enumerate}
\end{prop}
\begin{demo}
\begin{enumerate}
  \item  By Cauchy-Schwarz inequality as in \eqref{C.S}, for $f\in D(\Delta)$ we have
\begin{align*}
\mid(\tilde{\Delta} f,f)_{\beta^+}-(f,f)_{\beta^+}\mid&=\mid\sum_{x\in V}\sum_{y\in V}b(x,y)f(x)\overline{f(y)}\mid\\
&\leq\sum_{x\in V}\sum_{y\in V}b(x,y)\mid f(x)\overline{f(y)}\mid\\
&\leq(f,f)_{\beta^+}
\end{align*}
which implies that $W(\tilde{\Delta})\subset D(1,1)$.
  \item If $\displaystyle{\sum_{x\in V}\beta^{+}(x)}=\displaystyle{\sum_{x\in V}\beta^+(x)<\infty}$, the constant function is an eigenfunction of $\tilde{\Delta}$ associated to $0$. Then $0$ is an eigenvalue of $\tilde{\Delta}$. Now, we suppose that $f$ is an eigenfunction  of $\tilde{\Delta}$ associated to $0$, therefore $\big((\tilde{\Delta}+\tilde{\Delta}')f,f\big)=0$. Thus by connectedness of $G$, $f$ is constant.
\end{enumerate}
\end{demo}

It is obvious that $\mathcal{R}e(A)=\dfrac{1}{2}\big(A+A^{*}\big)$ if $A$ is a bounded operator, but this is not true in general. The result below establishes a link between the real part of a matrix and its eigenvalues considered as the roots of the characteristic polynomial, see \cite{L.G} page 8.
The adjoint of a square matrix is the transpose of its conjugate.
\begin{lemm}\label{lema}
Let $A$ be a square matrix of size $n$, $\lambda_k(A)$ and $\lambda_k(\mathcal{R}e(A))$, $k=1,..,n$ the eigenvalues of $A$ and $\mathcal{R}e(A)$ respectively. Suppose that the eigenvalues of $\mathcal{R}e(A)$ are labelled in the increasing order, so that, $\lambda_1(\mathcal{R}e(A))\leq\lambda_2(\mathcal{R}e(A))..\leq\lambda_n(\mathcal{R}e(A))$. Then
$$\sum^{n}_{k=n-q+1}\mathcal{R}e(\lambda_k(A))\leq \sum^{n}_{k=n-q+1}\lambda_k(\mathcal{R}e(A)),~~\forall q=1,..,n$$
and the equality prevails for $q=n$.
\end{lemm}
\begin{rem}
It should be noted that for a matrix $A$, 
$\lambda_k(\mathcal{R}e(A))$ and  $\mathcal{R}e(\lambda_k(A))$ are not equal in general.
We can see \cite{Ma} for a counter-example.
\end{rem}
In the following we study some generalities of eigenvalues of $\tilde{\Delta}^{D}_\Omega$, where $\Omega$ is a finite subset of $V$. We assume that they are ordered as follows:
$$\mathcal{R}e(\lambda_1(\tilde{\Delta}^{D}_\Omega))\leq\mathcal{R}e(\lambda_2(\tilde{\Delta}^{D}_\Omega))..\leq\mathcal{R}e(\lambda_n(\tilde{\Delta}^{D}_\Omega)).$$
%We prove the following  inequality involving $\lambda_1(\tilde{\Delta}^{D}_\Omega)$ and the lower bound of the spectrum of its real part.
\begin{lemm}\label{ag} Let $\Omega$ be a finite non-empty subset of $V$, we have
$$\lambda_1(\mathcal{R}e(\tilde{\Delta}^{D}_\Omega))\leq\mathcal{R}e(\lambda_1(\tilde{\Delta}^{D}_\Omega)).$$
\end{lemm}
\begin{demo}
Let $f$  be an eigenfunction associated to $\lambda_1(\tilde{\Delta}^D_\Omega)$.
By the variational principle of $\lambda_1(\tilde{H}^D_\Omega)$, we have
\begin{align*}
\lambda_1(\tilde{H}^D_\Omega)&\leq\dfrac{(\tilde{H}^D_\Omega f,f)_m}{(f,f)_m}\\
&=\dfrac{(\tilde{\Delta}^D_\Omega f,f)_m}{(f,f)_m}+\dfrac{\overline{(\tilde{\Delta}^D_\Omega f,f)}_m}{(f,f)_m}\\
&=\lambda_1(\tilde{\Delta}^D_\Omega)+\overline{\lambda_1(\tilde{\Delta}^D_\Omega)}.
\end{align*}
\end{demo}

 The next statement contains an additional information about the eigenvalues of $\tilde{\Delta}^D_\Omega$.
\begin{prop}\label{finit}
Let $\Omega$ be a finite non-empty subset of $V$ $(\#\Omega=n)$ such that $\partial_V\Omega\neq\emptyset$. Then the following assertions are true
\begin{enumerate}
\item $0<\mathcal{R}e(\lambda_1(\tilde{\Delta}^{D}_\Omega))\leq 1$.
\item $\lambda_1(\mathcal{R}e(\tilde{\Delta}^{D}_\Omega))+\lambda_n(\mathcal{R}e(\tilde{\Delta}^{D}_\Omega))\leq 2$.
\end{enumerate}
\end{prop}
\begin{demo}
\begin{enumerate}
\item From Theorem 4.3 of \cite{A.G}, we have $\lambda_1(\mathcal{R}e(\tilde{\Delta}^{D}_\Omega))>0$ and by Lemma \ref{ag} we conclude the left inequalty. Next, by Lemma \ref{lema} we have for $q=n$:
$$\sum^{n}_{k=1}\mathcal{R}e(\lambda_k(\tilde{\Delta}^{D}_\Omega))= \sum^{n}_{k=1}\lambda_k(\mathcal{R}e(\tilde{\Delta}^{D}_\Omega)$$
then
$$n\mathcal{R}e(\lambda_1(\tilde{\Delta}^{D}_\Omega))\leq\sum^{n}_{k=1}\lambda_k(\mathcal{R}e(\tilde{\Delta}^{D}_\Omega)=Tr\big(\mathcal{R}e(\tilde{\Delta}^{D}_\Omega)\big)=n$$
which proves that
$$\mathcal{R}e(\lambda_1(\tilde{\Delta}^{D}_\Omega))\leq 1.$$

\item It is deduced from the result of the symmetric case, see Theorem 4.3 \cite{A.G}.
\end{enumerate}
\end{demo}
\begin{coro} Let $\Omega$ be a finite non-empty subset of $V$, then
$$\mathcal{R}e(\lambda_n(\tilde{\Delta}^{D}_\Omega))<2.$$
\end{coro}
\begin{demo}
Applying the Lemma \ref{lema} for $q=1$, we get
$$\mathcal{R}e(\lambda_n(\tilde{\Delta}^{D}_\Omega))\leq\lambda_n(\mathcal{R}e(\tilde{\Delta}^{D}_\Omega)).$$
But by (2) of Proposition \ref{finit}, we have :
$$\lambda_n(\mathcal{R}e(\tilde{\Delta}^{D}_\Omega))\leq 2- \lambda_1(\mathcal{R}e(\tilde{\Delta}^{D}_\Omega)).$$
Then from the general property $\lambda_1(\mathcal{R}e(\tilde{\Delta}^{D}_\Omega))>0$, we conclude that $\lambda_n(\mathcal{R}e(\tilde{\Delta}^{D}_\Omega))< 2$.
\end{demo}
%\begin{rem}
%\begin{enumerate}
%\item $\Delta$ is not a normal operator in general.
%\item $\Delta$ is normal on the simple graph $\mathbb{Z}=\{..,-n,-n-1,..,0,..,n,n+1,..\}$ such that $\vec{E}=\{(k,k+1),~~ k\in \mathbb{Z}\}$, in fact:\\ $f\in D(\Delta)=D(\Delta^{*})=\ell^{2}(V,1)$ and we have for all $k\in \mathbb{Z}$, $$\Delta\Delta^{*}f(k)=\Delta^{*}\Delta f(k)=2f(k)-f(k-1)-f(k+1).$$
%\item  $\sigma_p(\Delta)=\sigma_p(\Delta^{*})$, in fact:\\Let $f$ be the eigenfunction associated to $\lambda$ such that $\Delta f=\lambda f$, so  $\Delta \overline{f}=\overline{\lambda} \overline{f}$. Hence, if $\lambda\in\sigma_p(\Delta)$ then $\overline{\lambda}\in\sigma_p(\Delta)$. Or it is known that $\sigma_p(\Delta^{*})=\{\overline{\lambda} ,~~\lambda\in \sigma_p(\Delta)\}$, where $\sigma_p(\Delta)$ is the set of eigenvalues of $\Delta$.
%\end{enumerate}
%\end{rem}

\section{Spectral study of the unbounded case}\label{ch}
This part includes the study of the  bounds on the numerical range and the essential spectrum of a closed Laplacian. Both issues can be approached via isoperimetric inequalities.\\

  %Often a problem translates by the data of an operator who is not closed. Rarely, we can apply explicitly certain properties for an unbounded non closed operator.
\subsection{Closable operator}
The purpose of the theory of unbounded operators is essentially to construct closed extensions of a given operator and to study their properties.
\begin{defi} \textbf{Closable operators}:
A linear operator $T : D(T) \to \mathcal{H}$ is closable if it has closed extensions.
\end{defi}
An interesting property for the Laplacian $\Delta$ is its closability.
\begin{prop}
Let $G$ be a graph satisfying the Assumption $(\beta)$. Then $\Delta$ is a closable operator.
\end{prop}
\begin{demo}
We shall use the Theorem of T. Kato which says that an operator densely defined is closable if its numerical range is not the whole complex plane, see \cite{Kts}, page 268.
Let $\lambda\in W(\Delta)$, there is $f\in \mathcal{C}_c(V)$ such that $\parallel f \parallel_m=1$ and $\lambda=(\Delta f,f)_m$. From the Green's formula we have, $$2 \mathcal{R}e(\lambda)=\sum_{(x,y)\in \vec{E}} b(x,y)\mid f(x)-f(y)\mid^{2}\geq 0.$$
It follows that $W(\Delta)\subset \big\{\lambda\in \mathbb{C},~~\mathcal{R}e(\lambda)\geq 0\big\}\subsetneq \mathbb{C}$.
\end{demo}

For such operators, another property of interest is the property of being closed.
\begin{defi}
The closure of $\Delta$ is the operator $\overline{\Delta}$, defined by
\begin{itemize}
\item $D(\overline{\Delta})=\big\{f\in \ell^{2}(V,m),~~\exists~ (f_n)_{n\in\mathbb{N}}\in \mathcal{C}_c(V),~f_n\to f \text{ and } ~\Delta f_n \text{ converge } \big\}$
\item $\overline{\Delta}f:=\displaystyle{\lim_{n\to \infty}}\Delta f_n ,~~ f\in D(\overline{\Delta})$ and $(f_n)_n\in \mathcal{C}_c(V)$ such that $f_n\to f$.
\end{itemize}
\end{defi}
For an unbounded operator the relation between the spectrum and the numerical range is more complicated. But for a closed operator we have the following inclusion, see \cite{Kts} and \cite{Arl}.
\begin{prop}\label{spn}
Let $T$ be a closed operator. Then $ \sigma_{ess}(T)\subset \overline{W(T)}$.
\end{prop}
 More precisely, let us define the following numbers:
 \begin{equation*}
\eta(T)=\inf\{\mathcal{R}e \lambda:~~\lambda \in \sigma(T)\}.
\end{equation*}
\begin{equation*}
\nu(T)=\inf\{\mathcal{R}e \lambda:~~\lambda\in W(T)\}.
\end{equation*}
\begin{equation*}
\eta^{ess}(T)=\inf\{\mathcal{R}e \lambda:~~\lambda \in \sigma_{ess}(T)\}.
\end{equation*}
The Proposition \ref{spn} induces this Corollary.
\begin{coro}
\begin{equation}\label{rang}
\eta^{ess}(\overline{\Delta})\geq \nu(\overline{\Delta}).
\end{equation}
\end{coro}
\begin{rem}
If $\Delta$ is self-adjoint, then $\eta(\Delta)=\nu(\Delta)$. But this is not the case in general. 
\end{rem}
\subsection{Cheeger inequalities}
For a non symmetric graph $G$, we prove bound estimates on the real part of the numerical range of $\Delta$ in terms of the Cheeger constant.   %Similar bounds are related with a fundamental inequality established in \cite{Klr}, \cite{Fuji} and \cite{A.G} for a simple symmetric graph. We give a lower bound of  $Q_{\Delta^{D}_\Omega}(f)$ for a subset $\Omega $ of $V$ and $f\in \mathcal{C}_c(\Omega)$ in terms of isoperimetric inequalities.
 We use this estimation to characterize the absence of the essential spectrum of $\overline{\Delta}$.\\

First, we recall the definitions of the Cheeger constants on $\Omega\subset V$:
 \begin{equation*}
h(\Omega)~ = \inf_{U\subset \Omega\atop finite}\frac{b(\partial_E U )}{m(U)}
\end{equation*}
and
\begin{equation*}
\tilde{h}(\Omega)~ = \inf_{U\subset \Omega\atop finite}\frac{b(\partial_E U )}{\beta^+(U)}
\end{equation*}
where for a subset $U$ of $V$, 
$$b(\partial_E U)=\sum_{(x,y)\in \partial_E U}b(x,y)$$
$$\beta^+(U)=\sum_{x\in U}\beta^+(x) \text{ and }m(U)=\sum_{x\in U}m(x).$$
We define in addition:
$$m_{\Omega}=\inf \left\{\frac{\beta^+(x)}{m(x)},~~x\in \Omega\right\}$$
$$M_{\Omega}=\sup \left\{\frac{\beta^+(x)}{m(x)},~~x\in \Omega\right\}.$$
%The weight $\beta$ on $x\in V$ is:
%$$\beta(x)=\beta^+(x)+\beta^-(x)=2\beta^+(x).$$

Cheeger's Theorems had appeared in many works on symmetric graphs. They give estimations of the bottom  of the spectrum of the Laplacian in terms of the Cheeger constant. The inequality \eqref{coi} controls the lower bound of the real part of $\lambda \in W(\Delta^{D}_\Omega)$. 
\begin{theo}\label{chee}
Let $\Omega\subset V$, the bottom of the real part of $W(\Delta^{D}_\Omega)$ satisfies the following inequalities:
\begin{equation}\label{coi}
\frac{{h}^{2}(\Omega)}{8}~\leq~M_\Omega\nu(\Delta^{D}_\Omega)~\leq ~\frac{1}{2}M_\Omega h(\Omega).
\end{equation}
\end{theo}
\begin{demo}
From the works of J. Dodziuk \cite{Dod} and A. Grigoryan \cite{A.G}, we can deduce the following bounds of the symmetric quadratic form $Q_{\Delta^{D}_\Omega}$ on $\mathcal{C}_c(\Omega)$, 
 $$\frac{{h}^{2}(\Omega)}{8}~\leq~M_\Omega\inf_{\|f\|_m=1}Q_{\Delta^{D}_\Omega}(f)~\leq ~\frac{1}{2} M_\Omega h(\Omega).$$
Then using the equality \eqref{symm} we conclude our estimation.
\end{demo}

We deduce in particular the following inequalities.
\begin{coro} Let $\Omega\subset V$, we have
\begin{equation*}\label{c}
\frac{{\tilde{h}}^{2}(\Omega)}{8}~\leq~\nu(\tilde{\Delta}^{D}_\Omega)~\leq ~\frac{1}{2}\tilde{h}(\Omega).
\end{equation*}
\end{coro}
\begin{prop} Let $\Omega\subset V$ and $g\in\mathcal{C}_c(\Omega),~\|g\|_m=1 $. Let $\lambda=(\Delta^{D}_\Omega g,g)_m\in W(\Delta^{D}_\Omega)$. Then
\begin{equation}\label{hhh}
   m_{\Omega}\frac{\mathcal{R}e(\tilde{\Delta}^{D}_\Omega g,g)_{\beta^+}}{( g,g)_{\beta^+}}\leq2\mathcal{R}e(\lambda)\leq M_{\Omega}\frac{\mathcal{R}e(\tilde{\Delta}^{D}_\Omega g,g)_{\beta^+}}{( g,g)_{\beta^+}}.
\end{equation}
\end{prop}
\begin{demo}
We have for all $x\in \Omega$
$$m_{\Omega}m(x)\leq \beta^+(x) \leq M_{\Omega}m(x)$$
therefore 
$$m_{\Omega}(g,g)_m\leq ( g,g)_{\beta^+} \leq M_{\Omega}(g,g)_m$$
which implies that:
 \begin{equation*}\label{Dir}
m_{\Omega}\frac{\mathcal{R}e(\tilde{\Delta}^{D}_\Omega g,g)_{\beta^+}}{( g,g)_{\beta^+}}\leq
\frac{Q_{\Delta^{D}_\Omega}(g)}{2(g,g)_m}\leq M_{\Omega}\frac{\mathcal{R}e(\tilde{\Delta}^{D}_\Omega g,g)_{\beta^+}}{( g,g)_{\beta^+}}
 \end{equation*}
 because $(\Delta^{D}_\Omega g,g)_m=(\tilde{\Delta}^{D}_\Omega g,g)_{\beta^+}$, for all $g\in\mathcal{C}_c(\Omega).$
\end{demo}
\begin{coro}  Let $\Omega\subset V$, we have
\begin{equation}\label{10}
m_\Omega\frac{\tilde{h}^{2}(\Omega)}{8}~\leq~\nu(\Delta^{D}_\Omega).
\end{equation}
\end{coro}
We can also estimate the real part of any element of the numerical range of $\Delta^{D}_\Omega$ in terms of
the isoperimetric constant $\tilde{h}$.
\begin{coro} \label{spe}
For all $\Omega\subset V$ and $\lambda\in W(\Delta^{D}_\Omega)$ we have
 \begin{equation}
    m_{\Omega}\big(2-\sqrt{4-\tilde{h}^2(\Omega)}\big)\leq 2
    \mathcal{R}e(\lambda)\leq M_{\Omega}\big(2+\sqrt{4-{\tilde{h}^2(\Omega)}}\big).
 \end{equation}
\end{coro}
\begin{demo}
We follow the same approach as Fujiwara in Proposition 1 \cite{Fuji}, and we apply it to the symmetric Laplacian $\tilde{H}^D_\Omega=\tilde{\Delta}^D_\Omega+\tilde{\Delta'}^D_\Omega$, we obtain, for all $g\in \mathcal{C}_c(\Omega)$
$$2-\sqrt{4-\tilde{h}^2(\Omega)}\leq\frac{2\mathcal{R}e(\tilde{\Delta}^{D}_\Omega g,g)_{\beta^+}}{( g,g)_{\beta^+}}\leq 2+\sqrt{4-{\tilde{h}^2(\Omega)}}.$$

Hence we obtain the result by a direct corollary of the inequality \eqref{hhh}.
\end{demo}
\subsection{Absence of essential spectrum from Cheeger constant }\label{s ess}
This subsection is devoted to the study of the essential spectrum relative to the geometry of the weighted graph. We  evaluate the interest of the study of the numerical range of non self-adjoint operators. Indeed, the knowledge of the numerical range of the Laplacian brings an essential information on its essential spectrum. \\

 We provide the Cheeger inequality at infinity on a filtration of graph $G$.
 \begin{defi}
A graph $H=(V_H,\vec{E}_H)$ is called a subgraph of $G=(V_G,\vec{E}_G)$ if $V_H\subset V_G$
and $\vec{E}_H=\big\{(x,y) ;~x,y \in V_H~~\big\}\cap \vec{E}_G $.
\end{defi}
\begin{defi}%\cite{Arv}
 A filtration of $G=(V,\vec{E})$ is a sequence of finite connected subgraphs $\{G_n=(V_n,\vec{E}_n),~~n\in \mathbb{N} \}$ such that $G_n\subset G_{n+1}$ and:
 $$\displaystyle\bigcup_{n\geq 1} V_n=V.$$
\end{defi}
Let $G$ be an infinite connected graph and $\{G_n,~n\in \mathbb{N}\}$ a filtration of $G$.
Let us denote
$$m_\infty=\lim_{n\rightarrow \infty}~ m_{V_n^c}$$
$$M_\infty=\lim_{n\rightarrow \infty}~ M_{V_n^c}$$
The Cheeger constant at infinity is defined by:
$$h_\infty=\lim_{n\rightarrow \infty}h(V^{c}_n).$$
\begin{rem}
These limits exist in $\mathbb{R}^+\cup \{\infty\}$ because $m_{V_n^c}$, $M_{V_n^c}$ and $h(V_n^c)$ are monotone sequences.
\end{rem}
\begin{rem}
The Cheeger constant at infinity $h_\infty$ is independent of the filtration. Indeed it can be defined, as in \cite{Fuji} and \cite{Klr}, by
$h_\infty =\displaystyle{\lim_{K\to G}h(K^c)}$, where $K$ runs over all finite subsets because the
graph is locally finite.
\end{rem}
\begin{defi}\label{ends}
$G$ is called with \textit{heavy ends} if $m_\infty=\infty$.
\end{defi}
\begin{lemm}\label{cc} For any subset $\Omega$ of $V$ such that $\Omega^c$ is finite, we have
$$\nu (\overline{\Delta}^D_\Omega)=\nu(\Delta^D_\Omega).$$
\end{lemm}
\begin{demo}
It is easy to see that 
$$b=\inf_{\lambda\in W(\overline{\Delta}^D_\Omega)}\mathcal{R}e(\lambda)\leq\inf_{\lambda\in W(\Delta^D_\Omega)}\mathcal{R}e(\lambda)=a.$$
Let $f\in D(\overline{\Delta}^D_\Omega)=\{f\in D(\overline{\Delta}),~f(x)=0,~\forall ~x\in \Omega^{c}\} $ such that $\parallel f\parallel_m=1$. Hence there is a sequence $(f_n)\in \mathcal{C}_c(V)=D(\Delta)$ which converges to $f$ and $(\Delta f_n)$ converges to $\overline{\Delta}f$. It follows that  $g_n=\mathbf{1}_\Omega f_n=0$ on $\Omega^c$ and
$(\Delta^{D}_\Omega g_n)$ converges to $\overline{\Delta}^D_\Omega f$. So
$$a\leq  \nu(\Delta^D_U)\leq \mathcal{R}e(\Delta^D_U g_n,g_n)_m \underset{n\to\infty}{\longrightarrow } \mathcal{R}e(\overline{\Delta}^D_U f,f)_m$$
then
$$a\leq b.$$

\end{demo}
\begin{theo}\label{spee}
The essential spectrum of $\overline{\Delta}$ satisfies:
\begin{equation*}
\frac{h^{2}_\infty }{8}~\leq~M_\infty\eta^{ess}(\overline{\Delta})
%\mathcal{R}e(\lambda_1^{ess}(\overline{\Delta})).
\end{equation*}
and
\begin{equation}\label{Ses}
m_\infty\frac{\tilde{h}^{2}_\infty }{8}~\leq~\eta^{ess}(\overline{\Delta}).
%\mathcal{R}e(\lambda_1^{ess}(\overline{\Delta})).
\end{equation}
\end{theo}
\begin{demo}
 Let $\{G_n,~n\in \mathbb{N}\}$ be a filtration of $G$, from the inequality \eqref{rang} we get,
$$\nu(\overline{\Delta}^D_{V^c_n}) \leq\eta^{ess}(\overline{\Delta}^D_{V^c_n}) .$$
From Theorem 5.35 of T. Kato page 244 \cite{Kts}, the essential spectrum is stable by a compact perturbation, we obtain
$$\sigma_{ess}(\overline{\Delta})=\sigma_{ess}(\overline{\Delta}^D_{V^c_n}).$$
Therefore
$$\nu(\Delta^D_{V^c_n})\leq \eta^{ess}(\overline{\Delta}),$$
we use Theorem \ref{chee} and the equality \eqref{10}, then we find the result by
taking the limit at $\infty$.
\end{demo}

The following Corollary follows from Theorem \ref{spee}. It gives an important characterization for the absence of the essential spectrum especially it includes the case of rapidly branching graphs.
\begin{coro}
The essential spectrum of $\overline{\Delta}$ on a heavy end graph $G$ with $\tilde{h}_\infty>0$ is empty.
\end{coro}
\begin{demo}
The emptiness of the essential spectrum for $\overline{\Delta}$ on a graph with \textit{heavy ends} is an immediate
Corollary of the inequality \eqref{Ses},
then if $m_\infty=\infty$ where $\tilde{h}_\infty> 0$, we have $\sigma_{ess}(\overline{\Delta})=\emptyset.$
\end{demo}

\thanks{
\textbf{Acknowledgments}: I take this opportunity to express my gratitude to
my thesis directors Colette Ann\'e and Nabila Torki-Hamza for all the fruitful discussions, helpful suggestions and their guidance during this work. This work was financially supported by the "PHC Utique" program of the French Ministry of Foreign Affairs and Ministry of higher education and research and the Tunisian Ministry of higher education and scientific research in the CMCU project number 13G1501 "Graphes, G\'eom\'etrie et th\'eorie Spectrale". Also I like to thank the Laboratory of Mathematics Jean Leray of Nantes (LMJL) and the research unity (UR/13ES47) of Faculty of Sciences of Bizerta (University of Carthage) for their financial and their continuous support.} I would like to thank  the anonymous referee for the careful reading of my paper and the valuable comments and suggestions.

\end{document}